\documentclass[a4paper,12pt,oneside]{article}
\usepackage{graphicx}
\usepackage[latin1]{inputenc}
\usepackage{graphicx}
\usepackage{amssymb}
\usepackage[english,french]{babel}

\begin{document}
  \begin{center}
    \huge \strut {Accurate a posteriori error evaluation in the reduced basis method} \par
  \end{center}
\vspace{1cm}

\begin{center}
 {Fabien Casenave}\\
{Universit\'e Paris-Est, CERMICS, \'{E}cole des Ponts ParisTech, 6 \& 8 av Blaise Pascal,
77455 Marne-la-Vall\'{e}e Cedex 2, FRANCE}\\
\end{center}
\vspace{1cm}

\section*{Abstract}
In the reduced basis method, the evaluation of the a posteriori estimator can become very sensitive to round-off errors.
In this note, the
origin of the loss of accuracy is revealed, and a solution to this problem is proposed and illustrated on a simple example.

\section{Introduction}

Consider the following discrete variational form, depending on a parameter $\mu\in\mathcal{P}$:
find $u_\mu$ in a Hilbert space $\mathcal{V}$ such that $\forall v\in \mathcal{V}$, $E_\mu:a_\mu(u_\mu,v)=b(v)$, where
$a_\mu$ is a bilinear form and $b$ a linear form.
In a many queries context, a quantity of interest of the solution $q(u_\mu)$ has to be computed
for many values of $\mu\in\mathcal{P}$. In this note, the following assumptions are made for simplicity, but the
conclusions are general:
(i) the variational formulation is coercive, (ii) $a_\mu$ has a so-called affine dependence on the parameter $\mu$ so that
$a_\mu=a_0+\alpha_1(\mu)a_1$,
(iii) the quantity of interest is the solution itself: $q(u_\mu)=u_\mu$.

Reduced basis (RB) strategies consist in replacing $E_\mu$ by an easily computable surrogate
$\hat{E}_\mu$, that is made precise below (see \cite{Machiels}). We denote $\hat{u}_\mu$ the solution of $\hat{E}_\mu$,
$N$ the size of the matrix involved in the resolution of $E_\mu$, and $\hat{N}$ the size of the matrix involved in
the resolution of $\hat{E}_\mu$. The RB method consists in two steps: (i) An offline stage, where a basis, whose vectors are
solutions of $E_\mu$ for well-chosen
values of the parameter $\mu$, is constructed using, e.g., a greedy algorithm on the parameter. During this stage,
$\hat{N}$ problems of size $N$ are
solved, and some quantities related to the solutions are stored. (ii) An online stage, where the precomputed quantities are
used to solve $\hat{E}_\mu$ for many values of $\mu$. In this stage, an a posteriori error estimator $\mathcal{E}(\mu)$
is also computed to check the quality of the approximation. This is called certification. The a posteriori error estimator
verifies $\|u_\mu-\hat{u}_\mu\|_{\mathcal{V}}\leqslant \mathcal{E}(\mu):=\beta_\mu^{-1}\|G_\mu \hat{u}_\mu\|_{\mathcal{V}}$,
where $\beta_\mu$ is the coercivity constant of $a_\mu$ (or a lower bound of it) and
$G_\mu$ is the unique affine application from $\mathcal{V}$ to $\mathcal{V}$ such that
$\forall (u,v)\in\mathcal{V}^2$, $(G_\mu u,v)_{\mathcal{V}}=a_\mu(u,v)-b(v)$. In this note, we consider different
ways to compute the same quantity $\mathcal{E}(\mu)$. We distinguish between formulae to compute $\mathcal{E}(\mu)$ by adding
an index to $\mathcal{E}(\mu)$. Thus, $\mathcal{E}_1(\mu):=\beta_\mu^{-1}\|G_\mu \hat{u}_\mu\|_{\mathcal{V}}$ is 
the first formula for the estimator, directly given by the definition. Since $\mathcal{E}_1(\mu)$ requires the computation
of a size $N$ scalar product, this formula is not compatible with the constraint that the computations in the online
stage should be of complexity independent of $N$.\\

Suppose that a reduced basis of size $\hat{N}$ has been computed
in the offline stage, namely a family $(\mu_i,u_i)_{i=1,...\hat{N}}$, where each of the $N$-dimensional vectors $u_i$ is
the solutions of $E_{\mu_i}$. The reduced problem is a Galerkin procedure on the $u_i$ basis: find
$\hat{u}_\mu\in \textnormal{Span}\{u_1, ..., u_{\hat{N}}\}$
such that $\forall j\in\{1,...,\hat{N}\}$, $\hat{E}_\mu:a_\mu(\hat{u}_\mu,u_j)=b(u_j)$.
Writing $\hat{u}_\mu=\sum_{i=1}^{\hat{N}}{\gamma_i(\mu)u_i}$, the reduced linear system is $\hat{A}_\mu\hat{U}_\mu=\hat{b}$,
where $(\hat{A}_\mu)_{i,j}=a_\mu(u_i,u_j)$, $(\hat{b})_j=b(u_j)$ and $(\hat{U}_\mu)_i=\gamma_i(\mu)$. Using the
affine parameter dependance, the matrix $(\hat{A}_\mu)_{i,j}=a_0(u_i,u_j)+\alpha_1(\mu) a_1(u_i,u_j)$ is built and solved in complexity
independent of $N$, provided that the quantities $a_0(u_i,u_j)$ and $a_1(u_i,u_j)$ have been precomputed in the
offline stage.

Consider the Riesz isomorphism $J$ from $\mathcal{V}'$ to $\mathcal{V}$ such that $\forall l\in\mathcal{V}'$, $\forall u\in\mathcal{V}$,
$\left(Jl,u\right)_{\mathcal{V}}=l(u)$. The operator $G_\mu$ inherits the affine
dependance of $a_\mu$ in $\mu$ since, $\forall u\in\mathcal{V}$,
\begin{equation}
G_\mu u=-Jb(\cdot)+Ja_0(u,\cdot)+\alpha_1(\mu)Ja_1(u,\cdot)=:G_{00}+G_0u+\alpha_1(\mu) G_1 u, 
\end{equation}
where $b(\cdot)$ and $a_k(u,\cdot)$, $k\in\{0,1\}$, are elements of $\mathcal{V}'$. The a posteriori error estimator
is then written in the following compact form:
\begin{equation}
\label{eq:est2}
\mathcal{E}_2(\mu):=\beta_\mu^{-1}\left(\delta^2 + 2s^t x_\mu + x_\mu^tSx_\mu
\right)^{\frac{1}{2}},
\end{equation}
where $\delta=\|G_{00}\|_{\mathcal{V}}$, $s_I=(G_{00},G_k u_i)_{\mathcal{V}}$, ${x_\mu}_I=
\alpha_k(\mu)\gamma_i(\mu)$, $S_{I,J}=(G_k u_i,G_l u_j)_{\mathcal{V}}$
(with $I$ and $J$ re-indexing respectively $(k,i)$ and $(l,j)$, $0\leqslant k,l\leqslant 1$, $1\leqslant i,j\leqslant\hat{N}$)
and $\alpha_0=1$.
Provided that $\delta\in\mathbb{R}$, $s\in\mathbb{R}^{2\hat{N}}$, and $S\in\mathbb{R}^{2\hat{N}\times 2\hat{N}}$
(which are independent of $\mu$) have been precomputed in the offline stage,
$\mathcal{E}_2(\mu)$ is computed in complexity independent of $N$. This is what is typically used in RB implementations.

\section{Round-off errors and certification}
\label{sec:machine_error}

Canuto, Tonn and Urban \cite{Canuto} identified  that the evaluation of $\mathcal{E}_2(\mu)$
suffers in practice from a loss of accuracy, which they attributed to the square root in \ref{eq:est2}.
Herein, we show more precisely that this loss of accuracy comes from round-off errors. Indeed, when substracting
two real numbers within floating point arithmetics,
the number of lost significant digits equals the number of common decimals between the two reals.
For simplicity, we neglect the round-off errors introduced when solving $E_\mu$ and $\hat{E}_\mu$, so that
the vectors of the reduced basis $u_i$ and the reduced solutions $\hat{u}_\mu$ are considered free of round-off errors.
Therefore, we only consider round-off errors in the evaluation of $\mathcal{E}_1(\mu)$ and $\mathcal{E}_2(\mu)$ due to the
summations.
We define the machine precision $\epsilon$ by the maximal floating point representation relative error of real numbers:
$\left|\frac{fl(x)-x}{x}\right|\leqslant\epsilon$. Under these hypotheses, the smallest possible values that can be
practically computed for $\mathcal{E}_1(\mu)$ and $\mathcal{E}_2(\mu)$ using floating point arithmetics with machine precison $\epsilon$
is bounded below by respectively $\frac{\delta}{\beta_\mu}\epsilon$ and $\frac{\delta}{\beta_\mu}\sqrt{\epsilon}$. This is
supported numerically (see section \ref{sec:application}).\\

This observation is of paramount importance since the certification of a RB procedure cannot be better
than these values. In a successful RB procedure, the value of the estimator gets smaller as the size $\hat{N}$ of the reduced
basis increases.
Enriching the basis with a new vector improves the quality of the approximation introduced by the method. As a result, there
exists $\hat{N}_0$ such that, $\forall\hat{N}\geqslant\hat{N}_0$, $\forall\mu\in\mathcal{P}$, $\mathcal{E}(\mu)\leqslant
\frac{\delta}{\beta_\mu}\sqrt{\epsilon}$. If $\hat{N}\geqslant\hat{N}_0$, $\mathcal{E}_2(\mu)$ is no longer suitable for computing
the a posteriori error estimator.

We notice that increasing the machine precision from $\epsilon$ to $\epsilon^2$
enables the accuracy of $\mathcal{E}_2(\mu)$ to reach the one of $\mathcal{E}_1(\mu)$. Thus, the use of quadruple precision is a first
solution, checked numerically in section \ref{sec:application}. This is however not pratical since current computer
architectures are optimized for double precision. Another solution is to develop an alternative algorithm for the evaluation
of the estimator that still achieves the machine precision. This is the purpose of the next
section.

\section{The new procedure for a posteriori error evaluation}
Consider that a reduced basis of size $\hat{N}$ has been constructed. Let us denote $d=1+3\hat{N}+2\hat{N}^2$.
For a given $\mu$ and $\hat{u}_\mu\in \textnormal{Span}\{u_1, ..., u_{\hat{N}}\}$, we define $X(\mu)\in\mathbb{R}^{d}$
as the vector with components $(1,{x_\mu}_I,{x_\mu}_I {x_\mu}_J)$,
with $1\leqslant I \leqslant J \leqslant 2\hat{N}$. Using the symmetry of the matrix $S$, we can write the right-hand side
of \ref{eq:est2} as a linear form in
$X(\mu)$: $\sum_{p=1}^{d}{q_p X_p(\mu)}$,
where $q_p$ is independent of $\mu$ and $X_p(\mu)$ is the p-th component of $X(\mu)$.

During the offline stage, we take $d$ values, possibly random, $\mu_r$, $r\in\{1,...,d\}$ of the parameter
$\mu$. Then, we compute the
vectors $X(\mu_r)$ and, using the accurate formula $\mathcal{E}_1$ for the estimator, the quantities
$V_r:=\left(\beta_{\mu_r}\mathcal{E}_1(\mu_r)\right)^2$. 
Finally, we define $T\in\mathbb{R}^{d\times d}$ as the matrix whose columns are formed by the vectors $X(\mu_r)$ and we assume
that $T$ is invertible, which was the case in our simulations.

In the online stage, suppose that we want to evaluate the estimator for the parameter value $\mu$. We compute the vector
$X(\mu)$ and solve the linear system $T\lambda(\mu)=X(\mu)$, for $\lambda(\mu)\in\mathbb{R}^d$. We then have
$X(\mu) = \sum_{r=1}^{d}{\lambda_r(\mu) X(\mu_r)}$ and
\begin{equation}
\sum_{p=1}^d q_p X_p(\mu) = 
\sum_{p,r=1}^d q_p \lambda_r(\mu) X_p(\mu_r) =
\sum_{r=1}^d{\lambda_r(\mu) \left(\beta_{\mu_r}\mathcal{E}_1(\mu_r)\right)^2} = \sum_{r=1}^d{\lambda_r(\mu) V_r}.
\end{equation}
This yields the new formula for computing the estimator,
\begin{equation}
\mathcal{E}_3(\mu):=\beta_{\mu}^{-1}\left(\sum_{r=1}^d{\lambda_r(\mu) V_r}\right)^{\frac{1}{2}}.
\end{equation}
Quite importantly, we notice
that the additional cost is such that the quantity $\mathcal{E}_3(\mu)$ is still computed in complexity independent of $N$.\\

The quantity $\left(\beta_\mu \mathcal{E}_2(\mu)\right)^2$ is a sum of terms, whose
first one is fixed and equals $\delta^2$. On the contrary, $\left(\beta_\mu \mathcal{E}_3(\mu)\right)^2$ is a sum whose
terms $V_r$ are updated each time a vector is added to the reduced basis. Since $V_r$ is computed for $\mu_r\in\mathcal{P}$,
$\left(\beta_\mu\mathcal{E}(\mu)\right)^2$ and $V_r$, (or, at least,
$\max_{\mu\in\mathcal{P}}{\left(\beta_\mu\mathcal{E}(\mu)\right)^2}$ and
$\max_{1\leqslant r\leqslant d}{V_r}$, which are the important quantities in a greedy selection and 
online certification) are of the same orders of magnitude.

\section{Numerical illustration}
\label{sec:application}
Consider the equation $-u''+\mu u=1$ on $]0,1[$ with $u(0)=u(1)=0$,
where $\mu\geqslant 1$ is the parameter. The analytic solution is
$u(x)=-\frac{1}{\mu}\left(\cosh\left(\sqrt{\mu}x\right)-1\right)+
\frac{\cosh\left(\sqrt{\mu}\right)-1}{\mu\sinh\left(\sqrt{\mu}\right)}\sinh\left(\sqrt{\mu}x\right)$.
The Lax-Milgram theory is valid, the coercivity constant is $1$ in the $H^1$-norm.
The estimator is given by $\mathcal{E}(\mu)=\|G_\mu \hat{u}_\mu \|_{H^1\left(]0,1[\right)}$.
Lagrange P$_1$ finite elements are used, with uniform mesh cells of length 0.005.
The RB method is carried-out until a reduced basis of size $6$ is constructed.
\begin{figure}[ht]
\centerline{
\includegraphics [width=8cm] {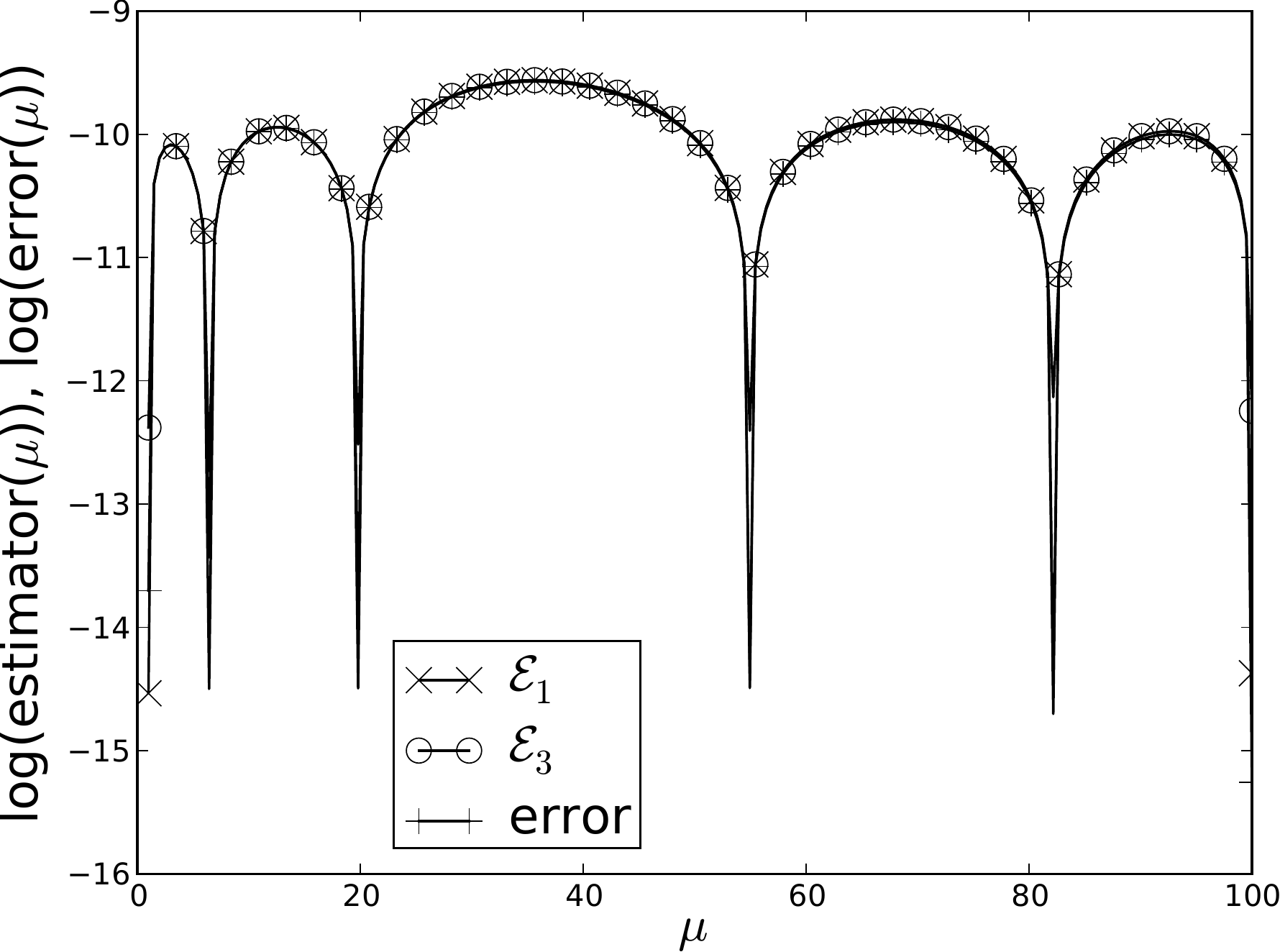}
\includegraphics [width=8cm] {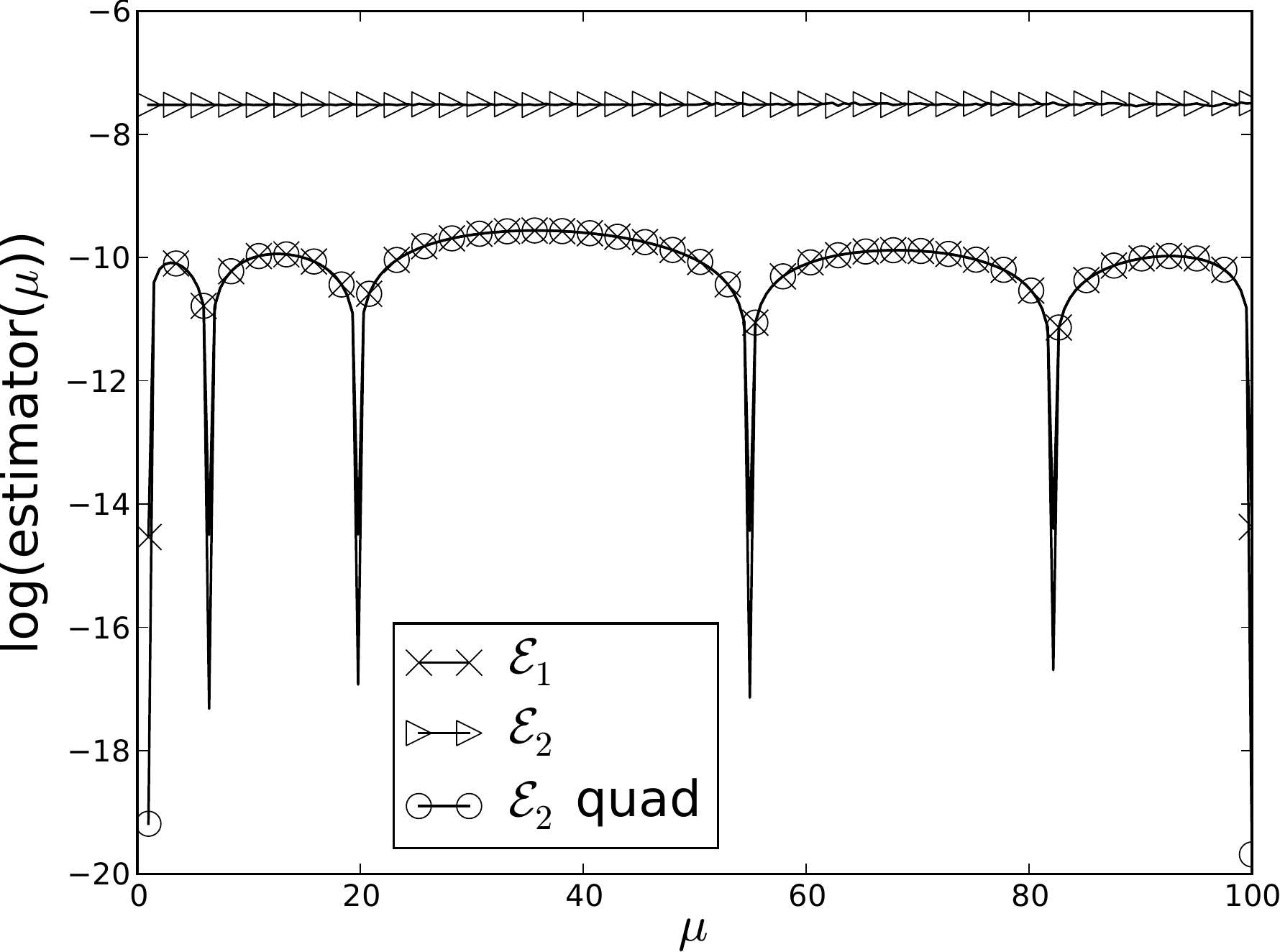}}
\caption{Left: $\mathcal{E}_1$, $\mathcal{E}_3$ algorithms for the estimator and error curves with
respect to the parameter $\mu$; Right: $\mathcal{E}_1$,
$\mathcal{E}_2$ (double and quadruple precision) curves.}
\label{fig:plot}
\end{figure}

On the left part of Figure \ref{fig:plot}, $\mathcal{E}_3$ yields the same curve as the accurate but expensive
$\mathcal{E}_1$ algorithm. Notice that the values of the a posteriori error estimators are very close to the values of the
error. This means that the efficiency of the estimator is very close to $1$.
On the right part of Figure \ref{fig:plot}, we see that $\mathcal{E}_2$ yields a flat curve for the estimator, meaning that
all information relative to the error is lost. As expected, the use of quadruple precision
enables $\mathcal{E}_2$ to recover the accuracy levels of $\mathcal{E}_1$. In this example,
$\frac{\epsilon\delta}{\beta}\approx 3\times 10^{-17}$ and $\frac{\sqrt{\epsilon}\delta}{\beta}\approx 3\times 10^{-9}$,
which sould be respectively compared to the numerical values of
$\mathcal{E}_2$ in quadruple precision and $\mathcal{E}_1$ ($10^{-17}$ and $10^{-15}$)
and $\mathcal{E}_2$ ($10^{-8}$).

\section{Conclusion}
To sum up, we have developed a procedure where the accuracy of the online
evaluation is limited by the accuracy of the evaluation of quantities precomputed during the offline stage, where heavy but
accurate algorithms are allowed. In the online stage, instead of a linear combination of $d$ terms, we have to
solve a linear system of size $d$, before doing a linear combination of the same size. We have increased
the accuracy of the estimator, with a procedure of complexity independent of the size $N$ of the initial problem.
When the size of the reduced basis increases, we observe that the condition number of the matrix $T$ increases as well.
Finally, we notice that oversampling strategies consisting in defining a least squares problem to compute $\lambda(\mu)$
such that $T$ is rectangular with more than $d$ columns
% (see oversampling methods in \cite{Champagnat})
improve the quality of our results when the RB is close to convergence. Experiments on a more complicated problem (external acoustics solved by integral equations where the criterion is on
the far field approximation of the diffracted acoustic potential) lead to similar conclusions, which will be
discussed in more detail elsewhere.

\section*{Acknowledgements}
This work was supported by EADS Innovation Works. The author thanks Nolwenn Balin, J\'{e}r\^{o}me Robert, Jayant Sen Gupta,
Guillaume Sylvand (EADS Innovation Works), and Alexandre Ern, Tony Leli\`evre (CERMICS) for fruitful discussions.


\begin{thebibliography}{ref}
\bibitem{Canuto}
{C. Canuto, T. Tonn and K. Urban, A posteriori error analysis of the reduced basis method for
nonaffine parametrized nonlinear pdes, SIAM J. Numer. Anal., 47:2001-2022 (2009)
}

\bibitem{Machiels}
{L. Machiels, Y. Maday, I.B. Oliveira, A.T. Patera and D.V. Rovas,
Output bounds for reduced-basis approximations of symmetric positive definite eigenvalue problems,
C. R. Acad. Sci. Paris, Ser. I 331 (2005)
}
% 
% \bibitem{Champagnat}
% {N. Champagnat, C. Chipot and E. Faou,
% Reconciling alternate methods for the determination of charge distributions: A probabilistic approach to high-dimensional least-squares approximations
% ArXiv e-prints, 1006.4996 (2010)
% }

\end{thebibliography}
\end{document}